\documentclass[11pt,a4paper,twoside]{amsart}
\usepackage{geometry}
\usepackage[utf8]{inputenc}
\usepackage[english]{babel}
\usepackage[T1]{fontenc}
\usepackage{lmodern}

\usepackage{amsmath}
\usepackage{amssymb}
\usepackage{amsthm}

\usepackage{mathtools}
\usepackage{stmaryrd}

\usepackage{wrapfig}
\usepackage{subcaption}

\usepackage{enumitem}
\usepackage{verbatim}
\usepackage{wrapfig}
\usepackage{makecell}
\usepackage{xcolor}

\usepackage{hyperref}
\hypersetup{colorlinks,linkcolor={blue!75!black},citecolor={blue!75!black},urlcolor={blue!75!black}}

\usepackage{tikz}

\DeclareMathOperator*{\supp}{supp}

\newcommand{\norm}[1]{\|#1\|}

\newcommand{\R}{\mathbb{R}}

\mathtoolsset{showonlyrefs}

\theoremstyle{theorem}

\newtheorem{theorem}{Theorem}[section] %

\newtheorem{lemma}[theorem]{Lemma}

\theoremstyle{definition}
\newtheorem{question}[theorem]{Question}

\theoremstyle{remark}
\newtheorem{remark}[theorem]{Remark}

\numberwithin{equation}{section}

\newcommand{\abs}[1]{\left\vert#1\right\vert}

\newcommand{\dint}{\,\mathrm{d}}

\title{Some new decay estimates for $(2+1)$-dimensional degenerate oscillatory integral operators}
\author{Yuxin Tan and Shaozhen Xu}

\address{School of Computer Science and Mathematics, Fujian University of Technology, Fuzhou 350118, China}
\email{tanyuxin14@mails.ucas.ac.cn}
\address{School of Information Engineering, Nanjing Xiaozhuang University, Nanjing 211171, China}
\email{shaozhen@nju.edu.cn}
\begin{document}
	
	\begin{abstract}
	In this paper, we consider the $(2+1)-$dimensional oscillatory integral operators with cubic homogeneous polynomial phases, which are degenerate in the sense of \cite{Tan06}. We improve the previously known $L^2\to L^2$ decay rate to $3/8$ and also establish a sharp $L^2\to L^6$ decay estimate based on fractional integration method.
	\end{abstract}
	
	\maketitle
	
	\setcounter{tocdepth}{1}
	
	\tableofcontents
	\section{Introduction}
	In this paper, we consider the following $(2+1)-$dimensional oscillatory integral operators
	\begin{equation}\label{Cubic-OIO}
		T_\lambda f(x,y)=\int_\R e^{i\lambda \left[P_1(x,y)t^2+P_2(x,y)t\right]}\psi(x,y,t)f(t)\dint t,
	\end{equation}
	where $\lambda>0$ and each $P_j(x,y)$ is a homogeneous polynomial of degree $j$ for $j=1,2$.  We seek to obtain decay or even sharp decay estimates in Lebesgue spaces for these operators. This goal is a small part of the project aiming to give sharp decay estimates for oscillatory integral operators with arbitrarily degenerate phases in arbitrary dimensions. The $(1+1)-$dimensional cases have been understood well, readers may refer to \cite{PhoSte92,See93,PhoSte94,PhoSte97,PhoSte98,PhoSteStu01,Ryc01,Yan04,Gre05,Xia17,SXY19} for more details. 
	The higher dimensional cases are much less understood, further background and progress can be found in the \emph{Introduction} of \cite{Xu23} and the references therein. \\
	
	In this paper, we focus on the operator of the form \eqref{Cubic-OIO}. In \cite{Tan06}, the author proved that if $P_2(x,y)$ has no multiple factors when it is factored into linear terms over $\mathbb{C}$, then it follows the (nearly) sharp decay estimates
	\begin{equation}
		\norm{T_\lambda f}_{L^2(\R^2)}\leq C\lambda^{-\frac{1}{2}}\log \lambda \norm{f}_{L^2(\R)}.
	\end{equation}
	Thus, the following question naturally arises:
	\begin{question}
		What are the decay or sharp decay estimates if $P_2(x,y)$ has multiple factors?
	\end{question}

	Before the attempt of answering this question, some assumptions should be introduced. Since $P_1(x,y)$ is a linear term and $P_2(x,y)$ is a quadratic form having multiple factors, then we may write 
	\begin{equation}
	P_1(x,y)=ax+by,\quad P_2(x,y)=(cx+dy)^2 \quad a,b,c,d\in \R.
	\end{equation}
	Next, apart from the trivial case $P_1(x,y)\equiv 0, P_2(x,y)\equiv 0$, we divide these polynomials into four cases according to coefficients.
	\begin{enumerate}
		\item [\emph{Case1}]: $P_1(x,y)\equiv 0$ and $P_2(x,y)\neq 0$,\\
		\item [\emph{Case2}]: $P_2(x,y)\equiv 0$ and $P_1(x,y)\neq 0$,\\
		\item [\emph{Case3}]: $P_1(x,y)\neq 0, P_2(x,y)\neq 0$ and $\left|\begin{matrix}
			a & b\\
			c & d
		\end{matrix}\right|=0$,\\
		\item [\emph{Case4}]: $P_1(x,y)\neq 0, P_2(x,y)\neq 0$ and $\left|\begin{matrix}
			a & b\\
			c & d
		\end{matrix}\right|\neq 0$.
	\end{enumerate}

	If the oscillatory integral operators of the form \eqref{Cubic-OIO} have phases in \emph{Case1, Case2} or \emph{Case3}, then we can reduce them to $(1+1)-$dimensional operators by a change of variables. Now, we take \emph{Case 3} for example. Without loss of generality, we may assume $a\neq 0$, since
	\[\left|\begin{matrix}
		a & b\\
		c & d
	\end{matrix}\right|=0,\]
	if we set 
	\[u=ax+by, v=y\]
	then 
	\[P_1(x,y)=u, P_2(x,y)=l^2u^2, \quad l=c/a.\]
	Recall the definition of $T_\lambda$, it can be  rewritten as
	\[\widetilde{T_\lambda} f(u,v)=\int_\R e^{i\lambda \left(ut^2+l^2u^2t\right)}\widetilde{\psi}(u,v,t)f(t)\dint t,\]
	here $\widetilde{T_\lambda}$ are essentially $(1+1)-$dimensional oscillatory integral operators. 
	From the result of Phong and Stein \cite{PhoSte94}, the sharp decay estimates follows
	\begin{equation}\label{CompEdgeEst}
		\norm{T_\lambda f(x,y)}_{L^2(\dint x\dint y)}=C_{a,b,c,d}\norm{\widetilde{T_\lambda}f(u,v)}_{L^2(\dint u\dint v)}\leq C_{\psi,a,b,c,d}\lambda^{-\frac{1}{3}}\norm{f}_{L^2(\R)}.
	\end{equation}
	Applying this argument to \emph{Case1} and \emph{Case2} shall lead to the following  sharp decay estimates
	\begin{equation}\label{VerEst}
		\norm{T_\lambda f(x,y)}_{L^2(\dint x\dint y)}\leq C_{\psi,a,b,c,d}\lambda^{-\frac{1}{4}}\norm{f}_{L^2(\R)}.
	\end{equation}
	
	Now we turn to the \emph{Case4}, in view of 
		\[\left|\begin{matrix}
		a & b\\
		c & d
	\end{matrix}\right|\neq 0,\]
	we set
	\[u=ax+by, v=cx+dy,\]
	then the original operator becomes
	\[\widetilde{T_\lambda} f(u,v)=\int_\R e^{i\lambda \left(ut^2+v^2t\right)}\widetilde{\psi}(u,v,t)f(t)\dint t.\] 
	It can be verified that
	\begin{equation}
		\norm{T_\lambda f(x,y)}_{L^2(\dint x\dint y)}=C_{a,b,c,d}\norm{\widetilde{T_\lambda}f(u,v)}_{L^2(\dint u\dint v).}
	\end{equation}
	To get the $L^2$ decay estimates for $T_\lambda$ in \emph{Case4}, it suffices to consider the operator of the form
	\begin{equation}
		\mathcal{O}_\lambda f(x,y)=\int_\R e^{i\lambda \left(xt^2+y^2t\right)}\psi(x,y,t)f(t)\dint t.
	\end{equation}

	In fact, for this operator, the second author \cite{Xu23} has given a $L^2\rightarrow L^2$ decay estimate,  with decay rate $1/4$, as a corollary of a sharp $L^4$ decay estimate which was established by employing broad-narrow method. In this paper, we will improve the $L^2\rightarrow L^2$ decay rate to $3/8$ and also establish a sharp $L^2\rightarrow L^6$ decay estimates based on fractional integration method. 
	
	The main theorem in this paper can be stated as follows.
	\begin{theorem}\label{MainThm}
		For the oscillatory integral operator $\mathcal{O}_\lambda$, we have the following decay estimates
		\begin{equation}\label{MainEstL2}
				\norm{\mathcal{O}_\lambda f}_{L^2(\R^2)}\leq C_{\psi}\lambda^{-\frac{3}{8}}\norm{f}_{L^2(\R)},
		\end{equation}
			\begin{equation}\label{MainEstL6}
		\norm{\mathcal{O}_\lambda f}_{L^6(\R^2)}\leq C_{\psi}\lambda^{-\frac{1}{4}}\norm{f}_{L^2(\R)},
	\end{equation}
and the second decay estimate is sharp.
	\end{theorem}

	\begin{remark}\label{Rem1}
		Following the example in \cite{Tan06}, one can show that the $L^2\to L^2$ decay rate should not be bigger than $1/2$, so we do not know whether \eqref{MainEstL2} is sharp or not. 
	\end{remark}
	\section{Improving $L^2\rightarrow L^2$ decay estimates}
	Before we give the proof of the inequality \eqref{MainEstL2}, we introduce a fundamental lemma due to Phong and Stein \cite{PhoSte94}, which plays a key role in the research of degenerate oscillatory integral operators. \\
	
	Consider the $(1+1)-$dimensional oscillatory integral operators of the form
	\[Tf(x)=\int_\R e^{i\lambda S(x,t)}\psi(x,t)f(t)dt\]
	where $S(x,t)$ is any polynomial and $\psi(x,t)$ is supported in a parallelogram with two sides parallel to the $t-$axis, having width $\delta$ in the $t-$direction. The following operator van der Corput Lemma gives us the locally uniform oscillatory estimates for $T$.
	\begin{lemma}\label{OscEst}
		Assume that
		\begin{enumerate}
			\item[\rm(i)] $\left|\partial_t^k\psi(x,t)\right|\leq B\delta^{-k}\left|\psi(x,t)\right|$, \quad $k=0,1,2$;\\
			\item[\rm(ii)] The degree of $S(x,t)$ in $t$ is less than or equal to $N$;\\
			\item[\rm(iii)] $0<\mu<\left|S_{xt}^{''}(x,t)\right|<A\mu.$
		\end{enumerate}
		Then it follows
		\begin{equation}
		\norm{T}_{L^2\to L^2}\leq C(\lambda\mu)^{-\frac{1}{2}},
		\end{equation}
	where the constant $C$ depends only on $A, B, N$.
	\end{lemma}
	Integration by parts twice in $t$ variable and Bernstein type inequality for derivatives of $S(x,t)$ will lead to the proof, details can be found in \cite{PhoSte94}.
	
	To balance the local oscillatory estimate, we also need the following size estimate.
	\begin{lemma}\label{SizEst}
		If $\psi(x,t)$ is supported in a region with $x-$cross length at most $\delta_1$, $t-$cross length at most $\delta_2$, then 
		\begin{equation}
			\norm{T}_{L^2\to L^2}\leq \norm{\psi}_{L^\infty} \left(\delta_1\delta_2\right)^{\frac{1}{2}}.
		\end{equation}
	\end{lemma}
This lemma is standard and can be seen a special case of Schur's test.
Now, we are ready to prove \eqref{MainEstL2}.
\begin{proof}[Proof of \eqref{MainEstL2}]
	To reduce the global estimate to a local one, we dyadically decompose $\mathcal{O}_\lambda$ into 
	\begin{align}
		\mathcal{O}_\lambda f(x,y)&=\sum_{j,k}\int_\R e^{i\lambda \left(xt^2+y^2t\right)}\psi(x,y,t)\phi_j(y)\phi_k(t)f(t)\dint t\\
		&:=\sum_{j,k}\mathcal{O}_{j,k}f(x,y)
	\end{align}
	where the inserted cut-off functions satisfy
	\[\sum_\sigma \phi_\sigma(\cdot)=1\]
	and $\supp \phi_\sigma\subset [2^{-\sigma-1},2^{-\sigma+1}], \sigma=j,k$.
	In the support of each $\mathcal{O}_{j,k}$, we know that
	\[\left|\partial_{xt}^{''}(xt^2+y^2t)\right|=|2t|\approx 2^{-k},\quad \left|\partial_{yt}^{''}(xt^2+y^2t)\right|=|2y|\approx 2^{-j}.\]
	Lemma \ref{OscEst} together with Lemma \ref{SizEst} implies
	\begin{align}
		\norm{\mathcal{O}_{j,k}}_{L^2\to L^2}&\lesssim \left(\lambda 2^{-k}\right)^{-\frac{1}{2}}\cdot 2^{-\frac{j}{2}},\\
		\norm{\mathcal{O}_{j,k}}_{L^2\to L^2}&\lesssim \left(\lambda 2^{-j}\right)^{-\frac{1}{2}},\\
		\norm{\mathcal{O}_{j,k}}_{L^2\to L^2}&\lesssim 2^{-\frac{j}{2}}\cdot 2^{-\frac{k}{2}}.	
	\end{align}
	Consequently,
	\begin{align}	
		\norm{\mathcal{O}_{\lambda}}_{L^2\to L^2}&\leq \sum_{j,k}\norm{\mathcal{O}_{j,k}}_{L^2\to L^2}\\
		&\lesssim \sum_{j,k}\min\left\{\left(\lambda 2^{-k}\right)^{-\frac{1}{2}}\cdot 2^{-\frac{j}{2}}, \left(\lambda 2^{-j}\right)^{-\frac{1}{2}}, 2^{-\frac{j}{2}}\cdot 2^{-\frac{k}{2}}\right\}\\
		&=\sum_{j>k/2} \min\left\{\left(\lambda 2^{-k}\right)^{-\frac{1}{2}}\cdot 2^{-\frac{j}{2}}, 2^{-\frac{j}{2}}\cdot 2^{-\frac{k}{2}} \right\}+\sum_{k\geq 2j}\min\left\{\left(\lambda 2^{-j}\right)^{-\frac{1}{2}}, 2^{-\frac{j}{2}}\cdot 2^{-\frac{k}{2}}\right\}\\
		&:=I+II.
	\end{align}
For $I$, we can see that
\begin{align}
	I&=\sum_{k=0}^\infty \min\left\{\left(\lambda 2^{-k}\right)^{-\frac{1}{2}}, 2^{-\frac{k}{2}}\right\}\sum_{j>k/2}^{\infty}2^{-\frac{j}{2}}\\
	&\approx \sum_{k=0}^\infty \min\left\{\left(\lambda 2^{-k}\right)^{-\frac{1}{2}}, 2^{-\frac{k}{2}}\right\}\cdot 2^{-\frac{k}{4}}\\
	&\approx \lambda^{-\frac{3}{8}}.
\end{align}
For $II$, we denote $k=2j+l$, so
\begin{align}
	II&=\sum_{l=0}^\infty\sum_{j=0}^\infty\min\left\{\left(\lambda 2^{-j}\right)^{-\frac{1}{2}}, 2^{-\frac{j}{2}}\cdot 2^{-j-\frac{l}{2}}\right\}\\
	&\approx \sum_{l=0}^{\infty}\lambda^{-\frac{3}{8}}2^{-\frac{l}{4}}\\
	&\approx \lambda^{-\frac{3}{8}}.
\end{align}
Thus we complete the proof.
\end{proof}
At the end of this section, we give an example to illustrate Remark \ref{Rem1}, in fact this kind of example has been given by \cite{PhoSte97} or in the same spirit in \cite{Tan06}. To make our argument be self-contained, we give the following relatively simple example.\\

Observe that the phase of $\mathcal{O}_\lambda$ is a homogeneous polynomial of degree 3, hence we may restrict each variable to a neighborhood of the origin with length less than $C\lambda^{-1/3}$ to make $\abs{\lambda xt^2+y^2t}\leq 1$. To complete the construction of the example, we further assume that the cut-off function satisfies
\begin{equation}\label{Supp-Func}
	\psi(x,y,t)=\begin{cases}
		0, &\quad \abs{(x,y,t)}\geq 1,\\
		1, &\quad \abs{(x,y,t)}\leq \frac{1}{2}.
	\end{cases}
\end{equation}
and the test function is
\begin{equation}\label{Tes-Func}
	f(t)=\chi_{[0,C\lambda^{-\frac{1}{3}}]}(t).
\end{equation}
If the $L^2\to L^2$ decay estimates with decay rate $\delta$ holds for any test function and any cut-off function, i.e.
\begin{equation*}
	\|\mathcal{O}_
	\lambda f\|_{L^2(\R^2)}\leq C_\psi\lambda^{-\delta} \|f\|_{L^2(\R)}.
\end{equation*}
then for the chosen test function and cut-off function, it obeys
\begin{equation*}
	\left[\iint\abs{\int_{\abs{t}\lesssim\lambda^{-1/3}}e^{i\lambda(xt^2+y^2t)}\psi(x,y,t)dt}^2dxdy\right]^{\frac{1}{2}}\lesssim C_\psi\lambda^{-\delta}\cdot\lambda^{-1/6}.
\end{equation*}
On account of the support function \eqref{Supp-Func}, we have
\begin{equation*}
	\lambda^{-\frac{2}{3}}\lesssim \left[\int_{\abs{y}\lesssim \lambda^{-1/3}}\int_{\abs{x}\lesssim \lambda^{-1/3}}\abs{\int_{\abs{t}\lesssim\lambda^{-1/3}}e^{i\lambda(x^2t+yt^2)}\psi(x,y,t)dt}^2dxdy\right]^{\frac{1}{2}}\leq C_\psi\lambda^{-\delta}\cdot\lambda^{-1/6}.
\end{equation*}
Since the inequality holds for arbitrarily large $\lambda$, then we have
\begin{equation*}
	\delta\leq \frac{1}{2}.
\end{equation*}

\section{Sharp $L^2\rightarrow L^6$ decay estimates}
This part is inspired by the endpoint Strichartz estimates of Keel and Tao \cite{KT98} and the proof of the 2-dimensional endpoint Stein-Tomas estimate using fractional integration method in \cite{MS13}. In continuity with the aim of \cite{Xu23}, we desire to explore the possibility of applying tools from Fourier restriction theory to degenerate oscillatory integral operators.\\

We briefly introduce the elements of our proof. The $TT^{*}$ method reduces the original oscillatory integral operators to operators with oscillatory kernels that can be bounded from above by the Stein-Wainger estimate \cite{SteWai01}. However, the dominant integral operator cannot be applied  the Hardy-Littlewood-Sobolev (HLS) inequality directly. We use a variant of the HLS inequality due to \cite{SY17} instead.
\begin{proof}[Proof of \eqref{MainEstL6}]
	We rewrite the estimate \eqref{MainEstL6} in a dual version
	\begin{equation}\label{DulMainEstL6}
				\norm{\mathcal{O}_\lambda^{*} g}_{L^2(\R)}\leq C_{\psi}\lambda^{-\frac{1}{4}}\norm{g}_{L^{\frac{6}{5}}(\R^2)}.
	\end{equation}
where the dual opeator $\mathcal{O}_\lambda^{*}$ is defined by
	\begin{equation}
		\mathcal{O}_\lambda^{*}g(t)=\int_{\R^2} e^{-i\lambda \left(xt^2+y^2t\right)}\bar{\psi}(x,y,t)g(x,y)\dint x\dint y.
	\end{equation}
Observe that
\begin{align}
	&\int_\R\left|\mathcal{O}_\lambda^{*}g(t)\right|^2\dint t\\
	&=\int_\R\mathcal{O}_\lambda^{*}g(t)\overline{\mathcal{O}_\lambda^{*}g(t)}\dint t\\
	&=\int_\R\left[\int_{\R^2} e^{-i\lambda \left(xt^2+y^2t\right)}\bar{\psi}(x,y,t)g(x,y)\dint x\dint y\right] \cdot \left[\int_{\R^2} e^{i\lambda \left(ut^2+v^2t\right)}\psi(u,v,t)\bar{g}(u,v)\dint u\dint v\right]\dint t\\
	&=\int_{\R^2}\int_{\R^2}\left[\int_\R e^{i\lambda\left[(u-x)t^2+\left(v^2-y^2\right)t\right]}\bar{\psi}(x,y,t)\psi(u,v,t)\dint t \right] g(x,y)\bar{g}(u,v)\dint u\dint v\dint x\dint y.
\end{align}
If we set
\begin{equation}
	K(u,v,x,y)=\int_\R e^{i\lambda\left[(u-x)t^2+\left(v^2-y^2\right)t\right]}\bar{\psi}(x,y,t)\psi(u,v,t)\dint t,
\end{equation}
then 
\begin{equation}
	\int_\R\left|\mathcal{O}_\lambda^{*}g(t)\right|^2\dint t=\int_{\R^2}\int_{\R^2}K(u,v,x,y) g(x,y)\bar{g}(u,v)\dint u\dint v\dint x\dint y.
\end{equation}
Write
\begin{equation}
	T_K\bar{g}(x,y)=\int_{\R^2}K(u,v,x,y) \bar{g}(u,v)\dint u\dint v.
\end{equation}
To prove \eqref{DulMainEstL6}, it suffices to verify that
\begin{equation}
	T_K:\quad L^{6/5}\rightarrow L^{6}.
\end{equation}
This reduction leads us to obtain a good estimate for the kernel which is a scalar oscillatory integral with a polynomial phase. This can be achieved using the following lemma due to Stein and Wainger \cite{SteWai01}.
\begin{lemma}
	Let $P(x)=\sum_{|\alpha|\leq D}c_\alpha x^\alpha$ be a $d-$dimensional polynomial of degree $\leq D$, $\varphi$ be a smooth function in the unit ball $B^d(0,1)$, and $\Omega$ be any convex subset of the unit ball. Then
	\begin{equation*}
		\left|\int_{\Omega}e^{iP(x)}\varphi(x)\dint x\right|\leq C_{D,d}\left(\sum_{0<|\alpha|\leq D}|c_{\alpha}|\right)^{-\frac{1}{D}}\sup_{x\in B^d(0,1)}\left(|\varphi(x)|+|\nabla\varphi(x)|\right).
	\end{equation*}
\end{lemma} 
This lemma implies 
\begin{equation}\label{KerEst}
	|K(u,v,x,y)|\leq C\lambda^{-\frac{1}{2}}\left(|u-x|+\left|v^2-y^2\right|\right)^{-\frac{1}{2}},
\end{equation}
where the constant is uniform with respect to all variables. Similar to what have been done in \cite{MS13}, we rewrite 
\begin{equation}
	T_K\bar{g}(x,y)=\int_\R T_{u,x}\bar{g_u}(y)\dint u,
\end{equation}
where the new operator inside the integral is
\begin{equation}
	T_{u,x}\bar{g_u}(y)=\int_{\R}K(u,v,x,y) \bar{g}(u,v)\dint v, 
\end{equation}
and $\bar{g}_u=\bar{g}(u,v)$.\\

Returning to the norm estimate of $T_Kg(x,y)$, we can see that
\begin{align}\label{UppBou}
	\norm{T_Kg(x,y)}_{L^6(\R^2)}^6&=\int_{\R}\int_{\R}\left|\int_\R T_{u,x}\bar{g}_u(y)\dint u\right|^6\dint x\dint y\\
	&\leq \int_{\R}\left[\int_{\R}\left(\int_\R\left| T_{u,x}\bar{g}_u(y)\right|^6\dint y\right)^{\frac{1}{6}}\dint u\right]^6\dint x.
\end{align}
The inequality is derived from Minkowski inequality. On account of \eqref{KerEst}, we can see that
\begin{equation}\label{InfToOne}
	\norm{T_{u,x}\bar{g}_u}_{L^\infty}\leq C\lambda^{-\frac{1}{2}}|u-x|^{-\frac{1}{2}}\norm{\bar{g}_u}_{L^1}.
\end{equation}
From \eqref{KerEst} we also have the pointwise estimate
\begin{equation}
	\left|	T_{u,x}\bar{g_u}(y)\right|\leq C\lambda^{-\frac{1}{2}}\int_{\R}\frac{|\bar{g}_u(v)|}{\left|v^2-y^2\right|^{\frac{1}{2}}}\dint v.
\end{equation}
Observe the right hand of the above inequality, it is in some sense similar to the usual fractional integration operators whereas the singularities are more sophisticated. In fact, this operator is a special case of the more general operators studied in \cite{SY17}. For convenience of the readers, we state the related mapping results given by Shi and Yan in \cite{SY17}.
\begin{lemma}\label{lemma-2}
	Let $W_{a,b}$ be the integral operator defined by
	\begin{equation*}
		W_{a,b}(f)(\eta)=\int_{\mathbb{R}^d}\Big||\xi|_a^a-|\eta|_a^a\Big|^{-\frac{d}{b}}f(\xi)d\xi.
	\end{equation*}For $d<a\leq b$,
	\begin{enumerate}
		\item $W_{a,b}$ is of weak type $(1,b/a)$.
		\item $W_{a,b}$ is bounded from $L^p(\mathbb{R}^d)$ to $L^q(\mathbb{R}^d)$, whenever $\frac{1}{p}=\frac{1}{q}+\frac{b-a}{b}$ and $1<p<\frac{b}{b-a}$.
	\end{enumerate}
\end{lemma}
This lemma implies the following $L^2\to L^2$ estimate
\begin{equation}\label{TwoToTwo}
	\norm{T_{u,x}\bar{g}_u}_{L^2}\leq C\lambda^{-\frac{1}{2}}\norm{\bar{g}_u}_{L^2}.
\end{equation}
Interpolation between \eqref{InfToOne} and \eqref{TwoToTwo} gives
\begin{equation}
	\norm{T_{u,x}\bar{g}_u}_{L^6}\leq C\lambda^{-\frac{1}{2}}|u-x|^{-\frac{1}{3}}\norm{\bar{g}_u}_{L^{6/5}}.
\end{equation}
Recall the upper bound of $\norm{T_Kg(x,y)}_{L^6(\R^2)}^6$ in \eqref{UppBou}, we know that
\begin{equation}
	\norm{T_Kg(x,y)}_{L^6(\R^2)}^6\leq C\lambda^{-3} \int_{\R}\left(\int_{\R}|u-x|^{-\frac{1}{3}}\norm{\bar{g}_u}_{L^{6/5}}\dint u\right)^6\dint x.
\end{equation}
Applying HLS inequality to the right hand of the above inequality gives
\begin{equation}
	\norm{T_Kg(x,y)}_{L^6(\R^2)}^6\leq C\lambda^{-3}\norm{g}_{L^{6/5}}^6.
\end{equation}
Consequently, H\"{o}lder's inequality implies
\begin{equation}
	\int_\R\left|\mathcal{O}_\lambda^{*}g(t)\right|^2\dint t\leq C\lambda^{-\frac{1}{2}}\norm{g}_{L^{6/5}}^2.
\end{equation}
This leads to the final result \eqref{MainEstL6}.\\

Next, we use an example in \cite{Xu23} to show the optimality of the decay estimate. In fact, the example is in spirit same with Knapp's example appearing in Fourier restriction theory. \\

Here we take \eqref{Supp-Func} as the cut-off function and choose the test function as
\begin{equation}\label{Tes-Func}
	f(t)=\chi_{[0,1]}(t).
\end{equation}
Since we assume the decay estimate
\begin{equation*}
	\|\mathcal{O}_
	\lambda f\|_{L^p(\R^2)}\leq C_\psi\lambda^{-\delta} \|f\|_{L^2(\R)}.
\end{equation*}
holds, then we know that
\begin{equation*}
	\left[\iint\abs{\int_0^1e^{i\lambda(xt^2+y^2t)}\psi(x,y,t)dt}^pdxdy\right]^{\frac{1}{p}}\leq C_\psi\lambda^{-\delta}.
\end{equation*}
Provided the support function \eqref{Supp-Func}, we have
\begin{equation*}
	\lambda^{-\frac{3}{2p}}\lesssim \left[\int_{\abs{y}\lesssim \lambda^{-1}}\int_{\abs{x}\lesssim \lambda^{-\frac{1}{2}}}\abs{\int_0^1e^{i\lambda(x^2t+yt^2)}\psi(x,y,t)dt}^pdxdy\right]^{\frac{1}{p}}\leq C_\psi\lambda^{-\delta}.
\end{equation*}
Since the inequality holds for arbitrarily large $\lambda$, then it requires
\begin{equation*}
	\delta\leq \frac{3}{2p}.
\end{equation*}
The case $p=6$ corresponds to our estimate, thus we complete the proof.
\end{proof}

	\section*{Acknowledgments}
	The first author is supported by Natural Science Foundation of
	Fujian Province (Grant No. 2021J011067) and  Foundation of Fujian University of Technology (Grant No. E0600421).
	The second author is supported by Jiangsu Natural Science Foundation under Grant No. BK20200308.


\end{document}